\documentclass[twocolumn,journal]{IEEEtran}
\usepackage[noadjust]{cite}

\ifCLASSINFOpdf
\usepackage[pdftex]{graphicx}
\else
\usepackage[dvips]{graphicx}
\fi

\usepackage{amsmath}
\usepackage{array}

\ifCLASSOPTIONcompsoc
  \usepackage[caption=false,font=normalsize,labelfont=sf,textfont=sf]{subfig}
\else
  \usepackage[caption=false,font=footnotesize]{subfig}
\fi
\usepackage{url}
\usepackage{amssymb}
\def\R{{\mathbb R}}
\DeclareMathOperator*{\argmin}{arg\, min}
\DeclareMathOperator{\sign}{sign}
\def\testmat{\text{te}}
\def\trainmat{\text{tr}}

\newtheorem{theorem}{Theorem}
\newtheorem{remark}{Remark}

\usepackage{algorithm}
\usepackage{algorithmicx}
\usepackage[noend]{algpseudocode}
\newcommand*\Let[2]{\State #1 $\gets$ #2}
\algrenewcommand\algorithmicrequire{\textbf{Input:}}
\algrenewcommand\algorithmicensure{\textbf{Output:}}

\usepackage{tikz}
\usepackage{color}
\definecolor{mygreen}{rgb}{0 0.7 0}
\definecolor{mygray}{rgb}{0.6 0.6 0.6}
\definecolor{myblue}{rgb}{0.7 0.7 1}
\definecolor{mypink}{rgb}{1 0 1}

\begin{document}

\title{Greedy Sensor Placement with Cost Constraints}

\author{Emily~Clark\thanks{E. Clark is with the Department of Physics at
the University of Washington, Seattle, WA, 98195-1560 USA email:
eclark7@uw.edu},
        Travis~Askham,
        Steven~L.~Brunton\thanks{S. Brunton is with the Department of Mechanical Engineering
at the University of Washington, Seattle, WA, 98195-2600 USA email:
sbrunton@uw.edu},~\IEEEmembership{Member,~IEEE,}
        J.~Nathan~Kutz\thanks{T. Askham and J. N. Kutz are with the Department of Applied
Mathematics at the University of Washington, Seattle, WA, 98195-3925 
USA email: askham@uw.edu and kutz@uw.edu},~\IEEEmembership{Member,~IEEE}}

\markboth{arXiv submission}%
{Clark \MakeLowercase{\textit{et al.}}: Sensor Placement with Cost Constraints}

\maketitle

\begin{abstract}
The problem of optimally placing sensors under a
cost constraint arises naturally in the design of 
industrial and commercial products, as well as in scientific 
experiments. We consider a relaxation of the full optimization 
formulation of this problem and then extend a well-established 
QR-based greedy algorithm for the optimal sensor placement problem 
without cost constraints. We demonstrate the effectiveness of 
this algorithm on data sets related to facial recognition, 
climate science, and fluid mechanics.
This algorithm is scalable and often identifies sparse 
sensors with near optimal reconstruction performance, while 
dramatically reducing the overall cost of the sensors.
We find that the 
cost-error landscape varies by application, with 
intuitive connections to the underlying physics.
Additionally, we include experiments for various pre-processing
techniques and find that a popular technique based on the
singular value decomposition is often sub-optimal.
\end{abstract}

\section{Introduction}

\IEEEPARstart{T}{he} problem of determining the optimal placement of sensors under
a cost constraint is relevant to many fields of scientific research and industry. 
Indeed, such considerations are critical in evaluating global monitoring systems and characterizing 
spatio-temporal dynamics (e.g. the brain, ocean and atmospheric dynamics, power grid networks, fluid flows, etc).  
For these applications, it is typical that only a limited number of measurements can be made of the system due to either prohibitive expense (i.e. either sensors are expensive, or they are expensive to place, or both) or the inability to place a sensor in a desired location (inaccessibility). 
Regardless, the goal of accurately reconstructing the state of the system from a limited number of measurements remains unchanged.  
To this end, we develop a principled, greedy sampling strategy whereby the sensor placement optimization is formulated as a cost-constrained problem in a relaxed form.  We further introduce a parameter representing the balance between the quality of the reconstruction and the cost, and thus can evaluate explicitly the cost-error curve.   The simple algorithmic structure proposed, which relies on a modification of the pivoted QR decomposition, provides an effective and scalable strategy for economical sensor placement for a wide range of scientific and engineering applications.


There is a significant body of literature on the subject
of signal reconstruction from a limited number of point
measurements (i.e. sensors). We do not seek to review this
literature here, but point to
\cite{wright2009robust,joshi2009sensor,krause2008near,
  krause2006near,zhang2008efficient,yildirim2009efficient,
  yang2010eof,manohar2017data} for a sense of the depth
and diversity of the existing mathematical optimization formulations. The related problem
of controlling a system based on a limited number of
measurements and actuators has also been well-studied. See, inter
alia, \cite{kammer1991sensor,martinez2006optimal,
lim1992method,cheng2008relay}.  
The critical extension of these techniques considered here involves the consideration
of the actual expense of the sensors and the cost of their placement, thus shaping a cost landscape that must be considered in order to more accurately assess the sensor placements.

In principle, the map from the measurements of a
system to the full-state reconstruction can take any
form. However, in this manuscript,
we view sensor placement as an interpolation problem,
i.e., given the values of a sample of some system at
its interpolation points (sensors), we would like to
approximately reconstruct the full state by applying a
{\em linear} map to these values. There are a few reasons
we limit ourselves to reconstruction based on a linear map:
(i) it is easy to check the stability and optimality of a
given set of points, (ii) it is straightforward to design
efficient algorithms for sensor placement, and, (iii)
as noted above, it is easy to interpret the sensor
locations as mathematical objects, i.e. they are
interpolation points.

Once the map is restricted to be linear,
the goal is then to find sensor locations which give
an accurate full-state approximation and result in a
stable interpolation map.
Of course, a brute-force solution 
of this problem may be obtained by searching over
all possible subsets of the sensors but this approach
quickly becomes intractable, as the number of subsets
increases combinatorially. We must then seek efficient
methods for finding nearly-optimal interpolation points.
The efficient computation of such points given 
samples of the system has a rich history and has
been considered in a variety of contexts. We will
quickly review some of the dominant themes of
the algorithms for selecting interpolation points.

Randomly placed sensors perform surprisingly
well. For instance, Wright et al. observed that,
given a generic basis in which samples of the signal will
be sparse, it is possible to reconstruct a signal
which has been downsampled or randomly projected
\cite{wright2009robust}. The compressed sensing
literature provides a theoretical basis for
the surprising effectiveness of random, or rather
{\em incoherent}, measurements in this setting; see,
inter alia, \cite{candes2006robust,donoho2006compressed,
  candes2006stable,candes2006near,baraniuk2007compressive}.
We note that such an approach does not necessarily
make use of any full-state observations of the system
(though some model for the system is implied) and
random sensors have been observed to be less efficient
than sensors which take this data into account
\cite{manohar2017data}.

A common data-driven approach is to start with
a {\em tailored basis} derived from the observed
samples, typically given by the dominant singular
vectors \cite{moore1981principal,willcox2002balanced,
  rowley2006model,juang1985eigensystem}.
See \cite{everson1995karhunen,willcox2006unsteady}
for early examples of signal reconstruction from
a limited number of sensors using such a basis.
Of course, random sensors may still be
used with tailored bases, but
better accuracy and stability are possible with sensors
chosen for the given basis. A number of heuristic
choices for the locations have been developed, including
placing sensors at the extrema of the singular
vectors \cite{zhang2008efficient,yildirim2009efficient,
  yang2010eof}.

Before we review some of the more principled
data-driven sensor selection algorithms, we require
some notation. For an
index set $J$ and any matrix ${\bf M}$, let ${\bf M}_{\cdot J}$
denote the matrix formed by the columns of ${\bf M}$
with index in $J$.
Let our samples of data be arranged in
the rows of a matrix ${\bf X}$ and let ${\bf \Psi}$ be
some matrix derived from ${\bf X}$ (e.g. ${\bf \Psi}$
may be taken to be the right singular
vectors of ${\bf X}$, random linear combinations of the
rows of ${\bf X}$, or ${\bf X}$ itself). If the size of
$J$ is fixed, it is known that
the set of indices $\hat{J}$ which maximizes the
product of the singular values of ${\bf \Psi}_{\cdot J}$
provides optimal interpolation points for ${\bf \Psi}$
\cite{gu1996efficient,cheng2005compression}
(see Theorem~\ref{thm:optimalid} below for the
definition of optimal).

The problem of finding
such a $\hat{J}$ is nonconvex and NP-hard, but there are
reasonable approximate algorithms.
Gu and Eisenstat developed a polynomial time algorithm
for computing $J$ when the optimality criterion is
relaxed slightly \cite{gu1996efficient}.
Joshi and Boyd formulated sensor placement as
an approximate convex problem, which may
be solved in polynomial time and is observed
to provide nearly
optimal sensors \cite{joshi2009sensor}.
While both of these approaches scale polynomially
in the number of sensors and the size of the
data, they are not as computationally efficient as some of the
existing greedy algorithms for interpolation, especially for high-dimensional data.
Further, the examples on which the greedy algorithms
are known to fail appear to be pathological, i.e.
it is incredibly unlikely that the greedy approach
will fail in practice.

The greedy sensor selection algorithm which is
of greatest interest in this manuscript is
based on the column pivoted QR decomposition.
In particular, for a given number of sensors $k$,
one simply selects $J$ to be the first $k$
column pivots of ${\bf \Psi}$ (see Section~\ref{sec:prelim} for
an explanation as to why this is a greedy approach
for maximizing the product of the singular
values of ${\bf \Psi}_{\cdot J}$). This algorithm is the
basis for practical approaches to computing the
interpolative decomposition \cite{cheng2005compression,
  martinsson2007interpolation}, which is commonly
used to compress low-rank matrices. The algorithm
is also used in the discrete empirical
interpolation method (DEIM) from reduced order
modeling \cite{barrault2004empirical,
  chaturantabut2010nonlinear},
in its more stable Q-DEIM formulation
\cite{drmac2016new}. For high dimensional problems
with many samples of data, standard techniques
from the burgeoning field of randomized algorithms
for linear algebra may be used to improve the
efficiency of these schemes
\cite{liberty2007randomized,halko2011finding}.

In the sensor placement techniques described above,
an optimal map and set of sensor locations are found
for a fixed number of sensors. This is equivalent to the
cost-constrained sensor placement problem when each
sensor has the same cost. In the
case that some sensor locations should be entirely
excluded, corresponding to an infinite cost, i.e. an inaccessible measurement location, and the
remaining locations are of uniform cost, again
the algorithms above may be used by simply narrowing
the search to the allowed sensor locations (note that
such a restriction has implications for the stability
of the interpolation map). 

We show that it is simple to modify the pivoted QR based
scheme to incorporate a cost constraint for problems
in between these extremes, i.e. for problems in which
some sensor locations cost more than others but may
be more informative. The method is obtained
by writing the cost-constrained problem in a
relaxed form, which introduces a parameter
representing the balance between the quality of
the reconstruction and the cost, and then
varying that parameter to trace out a cost-error
curve. For each value of the parameter, we use
a greedy algorithm to add sensor locations
one-by-one.

We test the performance of our methods on
data sets from facial recognition, climate
science, and fluid mechanics using a standard
training set/testing set apparatus. For some
data sets, the proposed
algorithm displays a significant advantage over
methods based on randomly selected sensors. We
also compare with known performance bounds and
brute-force answers when possible and find that
our algorithm is often near the optimal solution.

The remainder of this paper is organized as
follows. In Section~\ref{sec:prelim}, we summarize
some relevant results from the interpolation
literature and present our problem formulation.
We then develop an algorithm for sensor placement in
Section~\ref{sec:algorithm} which is a simple
extension of the existing methods. We include a brief
discussion of the effect of data pre-processing
(i.e. the choice of ${\bf \Psi}$) on the
quality of the sensor locations in Section~\ref{sec:projection};
in particular, we compare the performance when applied
to the raw data, the first several singular vectors
of the data, and randomized projections of the data.
We then apply our methods to three data sets
and discuss the performance in Section~\ref{sec:results}.
Finally, we provide some concluding remarks and indicate
possible future avenues for research in
Section~\ref{sec:conclusion}.



\section{Preliminaries and problem formulation}
\label{sec:prelim}

In this section, we fix some notation, formulate
the linear sensor placement problem with non-uniform
cost constraints, and review some of the literature on
the standard linear sensor placement problem with uniform cost.

\subsection{Setting, notation, and problem
  formulation}

Let ${\bf x}^i \in \R^n$ denote samples of some system
and let ${\boldsymbol \eta} \in \R_{+}^n$ denote non-negative 
costs associated with each sample location. We
collect the samples ${\bf x}^i$ as the rows of a matrix
${\bf X} \in \R^{m\times n}$. The sensor placement problem
with cost constraints then seeks out an optimal subset
$\hat{J}$ of the column indices of ${\bf X}$ which balances
the associated cost $\sum_{j\in \hat{J}} {\eta}_j$ with
reconstruction error and stability, which we define below.

For a given set of indices,
$J = \{j_1, \ldots, j_l \}$, it
is simple to construct the optimal linear map
for reconstructing the entries in ${\bf X}$. Let

\begin{equation}
  \hat{\bf T}(J) = \argmin_{{\bf T}\in\R^{l\times n}}
  \|{\bf X} - {\bf X}_{\cdot J} {\bf T} \|_F \; ,
\end{equation}
where ${\bf X}_{\cdot J}$ denotes the matrix given by
collecting the columns of ${\bf X}$ whose indices are
in $J$ and $\|\cdot \|_F$ denotes the Frobenius
norm. It is well known that
$\hat{\bf T}(J) = {\bf X}_{\cdot J}^\dagger {\bf X}$, where
${\bf X}_{\cdot J}^\dagger$ denotes the Moore-Penrose
pseudoinverse of ${\bf X}_{\cdot J}$, i.e. this is the least-squares solution. Therefore, the
relative reconstruction error for linear sensor
placement is given by

\begin{equation}
  e(J) = \frac{ \|{\bf X} - {\bf X}_{\cdot J} {\bf X}_{\cdot J}^\dagger {\bf X}\|_F}{\|{\bf X}\|_F}
\end{equation}
and the stability of the interpolation map
is determined by $\| {\bf X}_{\cdot J}^\dagger {\bf X} \|_\infty$.
In the following, we will then focus on
computing a subset $J$ such that the error
is small and the map is stable.

We will use much of the notation introduced
above throughout the paper. When necessary,
we will denote data matrices and errors
corresponding to a training set by
${\bf X}^\trainmat$ and $e^\trainmat(J)$ and
the matrices and errors corresponding to a
testing set by ${\bf X}^\testmat$ and $e^\testmat(J)$.
Note that

\begin{equation} \label{eq:errtest}
  e^\testmat(J) = \frac{ \|{\bf X}^\testmat - {\bf X}^\testmat_{\cdot J}
    {\bf X}_{\cdot J}^{\trainmat \dagger} {\bf X}^\trainmat\|_F}
  {\|{\bf X}^\testmat\|_F} \; ,
\end{equation}
i.e. the operator $\hat{\bf T}(J)$ is always
determined by the training set.
It is also common (see 
\cite{liberty2007randomized,halko2011finding,
chaturantabut2010nonlinear,manohar2017data}) 
to reduce the computational
cost associated with finding $J$ by applying
the algorithm to an $r\times n$ matrix
${\bf \Psi}^\trainmat$ which captures the dominant
features of ${\bf X}^\trainmat$ for some $r \ll m$,
e.g. to a matrix of singular vectors of
${\bf X}^\trainmat$ or a matrix given by random linear
combinations of the rows of ${\bf X}^\trainmat$.
In this case, we define
$\hat{\bf T}(J) = {\bf \Psi}_{\cdot J}^{\trainmat \dagger} {\bf \Psi}^{\trainmat}$.

It is now possible to define the
linear sensor placement problem with cost
constraints. Let $\boldsymbol{\eta}$ be the cost
vector as described above and let $s$ and
$b$ denote desired upper bounds on the
stability of the map $\hat{\bf T}(J)$ and the
budget, respectively. Then we may write
the cost-constrained problem as

\begin{equation} \label{eq:form1}
  \hat{J} = \argmin_{J} e(J) \text{ s.t. }
  \sum_{j\in J} {\eta}_j \leq b \text{ and }
  \|\hat{\bf T}(J)\|_{\infty,\text{vec}} \leq s \; ,
\end{equation}
where $\|\cdot\|_{\infty,\text{vec}}$ denotes the
maximum absolute value over the entries of
a matrix. Our algorithms will actually focus
on the following relaxation of \eqref{eq:form1}.
Note that, for a given $b$, there exists a
$\lambda$ such that the problem

\begin{equation} \label{eq:form2}
  \hat{J} = \argmin_{J} e(J) + \lambda
  \sum_{j\in J} {\eta}_j  \text{ s.t. }
  \|\hat{\bf T}(J)\|_{\infty,\text{vec}} \leq s \; 
\end{equation}
and \eqref{eq:form1} have the same solution.
Because we are often interested in the cost-error
landscape, we seek the solution of \eqref{eq:form2}
for a number of values of $\lambda$, so that we 
trace out a cost-error curve. We note that the 
main algorithm we present in 
Section~\ref{sec:algorithm} does not actually 
solve \eqref{eq:form2}. Instead, we seek a greedy 
approximate solution which does not strictly enforce 
the stability constraint but uses a heuristic strategy 
to bias the sensors in favor of stability.

In the case that all entries of $\boldsymbol{\eta}$ are equal
and positive, the constraint 
$\sum_{j\in J} {\eta}_j \leq b$ simplifies to the
constraint $|J| \leq b/{\eta}_1$. This is closely 
related to the problem of optimally placing a 
specified number of sensors, a well-studied problem
which we briefly review in the remainder of this 
section.

\subsection{Theoretical results}

As observed in \cite{cheng2005compression},
the strong rank revealing QR decomposition
methods of \cite{gu1996efficient} provide a
polynomial time algorithm for computing a
subset $J$ such that the error $e(J)$ is
near-optimally small and the map $\hat{\bf T}(J)$
is near-optimally stable. The optimal performance
is summarized in the following theorem.

\begin{theorem}{(adapted from \cite{cheng2005compression})}
\label{thm:optimalid}
  Let ${\bf X}\in \R^{m\times n}$ and let $k \leq l = \min (m,n)$.
  Then, there exists a $J$ such that $|J|=k$,

  \begin{equation}
    \|{\bf X} - {\bf X}_{\cdot J} \hat{\bf T}(J) \|_F \leq
    \sqrt{1 + k(l-k)} \sum_{j=k+1}^l \sigma_j({\bf X}) \; ,
  \end{equation}
  where $\sigma_j({\bf X})$ denotes the $j$th singular
  value of ${\bf X}$, and

  \begin{equation}
    \| \hat{\bf T}(J) \|_F \leq \sqrt{k(n-k)+k} \; .
  \end{equation}
\end{theorem}

Remarkably, if we relax the bounds to

\begin{equation}
  \|{\bf X} - {\bf X}_{\cdot J} \hat{\bf T}(J) \|_F \leq
  \sqrt{1 + lk(l-k)} \sum_{j=k+1}^l \sigma_j({\bf X}) \; ,
\end{equation}
and

\begin{equation}
  \| \hat{\bf T}(J) \|_F \leq \sqrt{nk(n-k)+k} \; ,
\end{equation}
then there exist algorithms which compute such
a $J$ using, typically, $O(mnk)$ flops and
at most $O(mnl)$ flops \cite{gu1996efficient}.
However, the standard QR algorithm with column
pivoting tends to achieve similar bounds, with
the worst case behavior limited to pathological
examples (see the Kahan matrix example of
\cite{gu1996efficient}).

\subsection{QR with column pivoting for sensor placement}
\label{sec:qrpivot}
Many of the standard implementations of routines
for computing a QR decomposition are based on the 
application of Householder reflectors to triangularize
a matrix \cite{businger1965linear}. Let a vector
${\bf v} \in \R^m$ be given. We can then define a Householder 
reflector \cite{householder1958unitary} which maps 
${\bf v}$ to $\|{\bf v}\|_2 {\bf e}^1$, where ${\bf e}^1$ is the first 
standard basis vector in $\R^m$. Let $\sigma = \|{\bf v}\|_2$
and let $v_1$ denote the first entry of ${\bf v}$. Then,
the matrix 

\begin{equation}
 {\bf H}({\bf v}) := {\bf I} - \frac{({\bf v}+\sign(v_1) \sigma {\bf e}^1)
({\bf v}+ \sign(v_1) \sigma {\bf e}^1)^\intercal}{\sigma(\sigma+|v_1|)}
\end{equation}
maps ${\bf v}$ to $-\sign(v_1)\sigma {\bf e}^1$. It is simple to 
verify that there exists a ${\bf u}\in \R^m$ of unit norm 
such that the map is given by ${\bf I}-2{\bf u}{\bf u}^\intercal$, which 
is the standard form of a Householder reflector. 
This expression also makes it clear that ${\bf H}({\bf v})$ is 
its own inverse.

Using the notation of \cite{gu1996efficient},
the standard QR decomposition algorithm with 
column pivoting applied to a matrix ${\bf X}\in \R^{m\times n}$
proceeds by iteratively defining an orthogonal 
matrix ${\bf Q}^k$ and a permutation of the indices
$J^k = \{j^k_1, \ldots, j^k_n\}$ such that 

\begin{equation} \label{eq:partialqr}
{\bf X} {\bf P}(J^k) =  {\bf Q}^k \begin{pmatrix} {\bf A}^k & {\bf B}^k \\
    0 & {\bf C}^k
  \end{pmatrix} \; ,
\end{equation}
where ${\bf A}^k \in \R^{k\times k}$ is upper triangular, 
${\bf B}^k \in \R^{k\times (n-k)}$, ${\bf C}^k \in \R^{(m-k)\times (n-k)}$,
and

\begin{equation}
 {\bf P}(J^k) = \begin{pmatrix} {\bf e}^{j^k_1} & {\bf e}^{j^k_2} & 
    \cdots & {\bf e}^{j^k_n}
  \end{pmatrix} \; .
\end{equation}

Let $J^0 = \{1,\ldots n\}$ and 
${\bf Q}^0 = {\bf I}$. Denote the columns of ${\bf C}^k$ by 
${\bf c}^{k,i}$ for $i = 1,\ldots, n-k$. To obtain the 
$k+1$st iterate from the $k$th, let ${\bf c}^{k,l}$ be
the column of ${\bf C}^k$ with the largest norm. 
We update the permutation indices to move this
column to the front, i.e. 
$j^{k+1}_{k+1} = j^k_{k+l}$, $j^{k+1}_{k+l} = j^k_{k+1}$,
and $j^{k+1}_i = j^k_i$ for all other $i$. If we 
update the orthogonal matrix ${\bf Q}^k$ via

\begin{equation}
  {\bf Q}^{k+1} = {\bf Q}^k \begin{pmatrix} {\bf I} & 0 \\ 0 & {\bf H}({\bf c}^{k,l})
  \end{pmatrix} \; ,
\end{equation}
then 

\begin{equation} \label{eq:partialqr}
{\bf X} {\bf P}(J^{k+1}) =   {\bf Q}^k \begin{pmatrix} {\bf A}^{k+1} & {\bf B}^{k+1} \\
    0 & {\bf C}^{k+1}
  \end{pmatrix} \; ,
\end{equation}
where ${\bf A}^{k+1}$, ${\bf B}^{k+1}$ and ${\bf C}^{k+1}$ are of the
correct form and ${\bf A}^{k+1}_{k+1,k+1} = \pm \|{\bf c}^{k,l}\|$.

The idea of using the column pivoted QR algorithm
for sensor placement is that, for a given $k$,
the first $k$ pivots $j^k_1, \ldots, j^k_k$ should
be a good choice of sensor locations. Let 
$J = \{ j^k_1, \ldots , j^k_k \}$. Then 

\begin{equation}
\hat{\bf T}(J) {\bf P}(J^k) = 
\begin{pmatrix} {\bf I} & 
\left( {\bf A}^k \right )^{-1} {\bf B}^k
\end{pmatrix}
\end{equation}
and $e(J) = \|{\bf C}^k \|_F / \|{\bf X}\|_F$, so that both 
the stability of the map $\hat{\bf T}(J)$ and the 
error $e(J)$ are determined by the factors ${\bf A}^k$, ${\bf B}^k$,
and ${\bf C}^k$. 

Of course, for any permutation $J^k$ of the indices,
it is possible to construct a ${\bf Q}^k$ using Householder
reflectors such that ${\bf A}^k$, ${\bf B}^k$, and ${\bf C}^k$ are of 
the correct form. In \cite{gu1996efficient},
it was shown that perfectly stable and highly 
accurate sensors, in the sense of Theorem~\ref{thm:optimalid},
may be obtained by maximizing $|\det {\bf A}^k|$ over all 
possible permutations $J^k$. Because the ${\bf A}^k$ are 
upper triangular, the column-pivoted QR procedure 
may be seen as a greedy method which approximates the
optimal solution iteratively, updating $J^{k+1}$ 
so that $| \det {\bf  A}^{k+1}|$ is as large as
possible with the first $k$ entries of 
$J^{k+1}$ fixed to be equal to the first $k$ 
entries of $J^k$.

\begin{remark}
  For the calculations in Sections~\ref{sec:projection}
  and \ref{sec:results},
we computed $\hat{\bf T}(J)$ using the formula 
$\hat{\bf T}(J) = {\bf X}_{\cdot J}^\dagger {\bf X}$ and found that
this worked well for our examples, though there may
be an advantage to using a more numerically stable
definition for $\hat{\bf T}(J)$, see (3.13) of 
\cite{cheng2005compression}.
\end{remark}

\subsection{Other work}

The reduced order modeling community has long 
used a different method for selecting interpolation
points of snapshot matrices: the discrete empirical
interpolation method (DEIM) 
\cite{chaturantabut2010nonlinear}. Let ${\bf v}^1,\ldots, {\bf v}^k$
denote the first $k$ right singular vectors of ${\bf X}$.
The original DEIM procedure sets the interpolation
points to be the pivots used when applying 
Gaussian elimination with partial pivoting to 
the matrix $({\bf v}^1, \ldots, {\bf v}^k)$. In \cite{drmac2016new},
the Q-DEIM procedure was defined which instead chooses 
the interpolation points as the first $k$ pivots of 
column-pivoted QR applied to the matrix 
${\bf \Psi} = ({\bf v}^1,\ldots,{\bf v}^k)^\intercal$. Q-DEIM has been
observed to be generally more stable and accurate than
DEIM \cite{drmac2016new,manohar2017data}. In this
manuscript, we will refer to the methods which define 
${\bf \Psi}$ in terms of singular vectors as DEIM-based.

For problems with many samples of data (large $m$),
methods from randomized linear algebra may be used 
to decrease the computational cost of sensor placement
\cite{liberty2007randomized,halko2011finding}. The
randomized linear algebra approach is based on the
fact that, with $r = k+t$ and the entries of
${\bf G}\in \R^{r\times m}$ drawn independently from a standard
normal distribution, the row space of the
matrix ${\bf \Psi} = {\bf G}{\bf X}$ closely approximates the space
spanned by the first $k$ right singular vectors of
${\bf X}$ with a failure probability that decreases
super-exponentially in $t$. For a matrix with
fast singular value decay, typically $t$ is taken
to be 10 \cite{halko2011finding}. For data matrices, 
the singular values often decay slowly, so we take 
$r =2k$ in our examples. We explore the performance 
of using true singular vectors, as in the DEIM 
framework, and using a randomized approach in
Section~\ref{sec:projection}.

\section{Algorithm for sensor placement
  under cost constraints} \label{sec:algorithm}

In this section, we present pseudo-code for a greedy
approach to the relaxed version of the cost-constrained
sensor placement problem \eqref{eq:form2}. The algorithm
is based on the column-pivoted QR algorithm described 
in Section~\ref{sec:qrpivot}, where the
pivot column is chosen to balance the decrease in the 
error $e(J)$ with the increase in the total cost
$\sum_{j\in J} {\eta}_j$. Within the pseudo-code, this 
balance is represented by the parameter $\gamma$ 
(see Algorithm~\ref{alg:qrpc}).

\begin{algorithm}[t]
  \caption[Caption]{\label{alg:qrpc} QR pivoting with cost constraints.}
  \begin{algorithmic}[1]
    \Require{${\bf X}$, $k$, $\boldsymbol{\eta}$, $\gamma$}
    \Let{$m,n$}{\texttt{size}$({\bf X})$}
    \Let{${\bf R}$}{\texttt{copy}$({\bf X})$}
    \Let{$\boldsymbol{\tau}$}{\texttt{zeros}$(\min (m,n),1)$}
    \Let{$J$}{$1:n$}
    \For{$i = 1,\ldots,k$}
    \For{$p = i,\ldots,n$}
    \Let{$\nu_p$}{$\|{\bf R}_{i:m,p}\|_2 
      - \gamma {\eta}_{j_p\vphantom{\hat{j_p}}}$}
    \EndFor
    \Let{$l$}{index of the maximum of $\boldsymbol{\nu}_{i:n}$}
    \Let{${\bf v}$}{${\bf R}_{i:m,i-1+l}$}
    \State \texttt{swap}$({\bf R}_{i:m,i},{\bf R}_{i:m,i-1+l})$
    \State \texttt{swap}$(j_i,j_{i-1+l})$
    \Let{$\sigma$}{$\|{\bf v}\|_2$}
    \Let{${\bf w}$}{$({\bf v}+\sign(v_1)\sigma {\bf e}^1)/\sqrt{\sigma (\sigma + |v_1|)}$}
    \Let{${\bf R}_{i:m,i:n}$}{${\bf R}_{i:m,i:n}-{\bf w}{\bf w}^\intercal {\bf R}_{i:m,i:n}$}
    \Let{$\tau_i$}{$w_1$}
    \Let{${\bf R}_{i+1:m,i}$}{${\bf w}_{2:\text{end}}$}
    \EndFor
    \Ensure{${\bf R}$, $\boldsymbol{\tau}$, $J$}
  \end{algorithmic}
\end{algorithm}

We note that, after iteration $i$ of the outer for-loop
in Algorithm~\ref{alg:qrpc}, 
the matrix ${\bf A}^i$ is stored in ${\bf R}_{1:i,1:i}$, ${\bf B}^i$ 
is stored in ${\bf R}_{1:i,i+1:n}$, and ${\bf C}^i$ is stored
in ${\bf R}_{i+1:m,i+1:n}$. Therefore, when $\gamma = 0$, 
the pivot chosen in Algorithm~\ref{alg:qrpc}
is the same as the pivot in column pivoted QR. 
At step $i+1$ of the outer for-loop, the difference
in the error, $e(J_{1:i})-e(J_{1:i+1})$, is at least 
$\|{\bf v}\|_2/(2\|{\bf X}\|_F\sqrt{n})$ (where ${\bf v}$ is as in the 
pseudo-code). Therefore,
a positive $\gamma$ balances the decrease 
in the error, $e(J_{1:i})-e(J_{1:i+1})$, with the cost 
of the pivot, ${\eta}_{j_{i+1}\vphantom{\hat{j_{i+1}}}}$. 

\begin{figure*}[t]
	\centering
	\includegraphics[width = 0.9\textwidth]{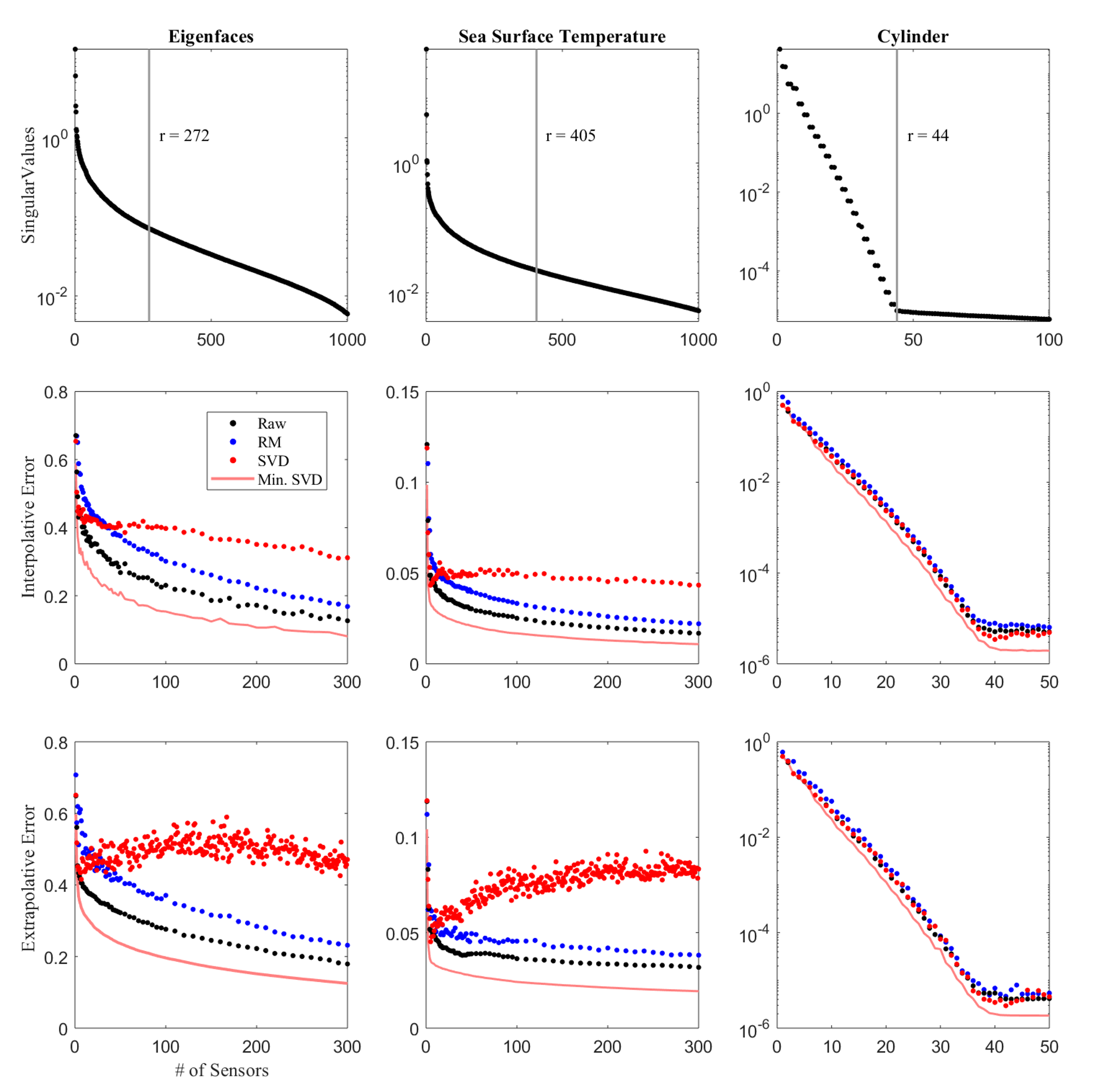}
		\vspace*{-.1in}
	\caption{Reconstruction error versus the number of sensors for the three data sets described in Section \ref{sec:results}, using $p$ SVD modes, random linear combinations with $2p$ modes (abbreviated RM in the legend), and raw data. The top row shows a log plot of the normalized singular value spectrum, with the vertical gray line indicating the Gavish-Donoho cutoff $r$ \cite{gavish2014optimal}. The remaining plots show the average reconstruction error given sensors placed using the three pre-processing methods discussed in the text. The first column provides eigenface results, the second gives sea surface temperature reconstruction errors, and the third shows reconstruction errors for the fluid flow behind a cylinder on a log scale. The middle row of the figure shows interpolative error, where the training set consists of a randomly-chosen subset of the data, while the bottom row gives the extrapolative data, which takes the first 80\% of the parameter space. All plots also show a rough estimate of the minimum error at a given number of sensors (the solid red line), obtained by projecting the full image onto the SVD modes.}
	\label{fig:Interp_Extrap}
\end{figure*}

\begin{remark}
  In Algorithm~\ref{alg:qrpc}, the pivot chosen at 
  step $i+1$ of the outer for-loop does not necessarily 
  correspond to the natural greedy choice, i.e. the pivot
  which minimizes $e(J_{1:i+1}) + \gamma \sum_{j\in J_{1:i+1}} 
  {\eta}_j$ with $J_{1:i}$ fixed. Such a pivot 
  could be computed, though at greater cost than 
  the present algorithm. Further, there  
  is another reason to avoid such a strategy: 
  it completely ignores the stability of 
  the resulting map. By instead pivoting based on 
  column size, we bias the algorithm toward choosing
  stable pivots while still incorporating some sense 
  of the reduction in error. 
\end{remark}

\begin{remark}
  We have also implemented an analogous algorithm
  based on Gaussian elimination with partial 
  pivoting (in the spirit of the original DEIM algorithm). 
  Because this approach performs significantly
  worse than the QR-based algorithm, we omit the
  details.
\end{remark}

\section{Data, singular vectors, and random projections} \label{sec:projection}

\begin{figure*}[t]
	\centering
	\includegraphics[width = 0.9\textwidth]{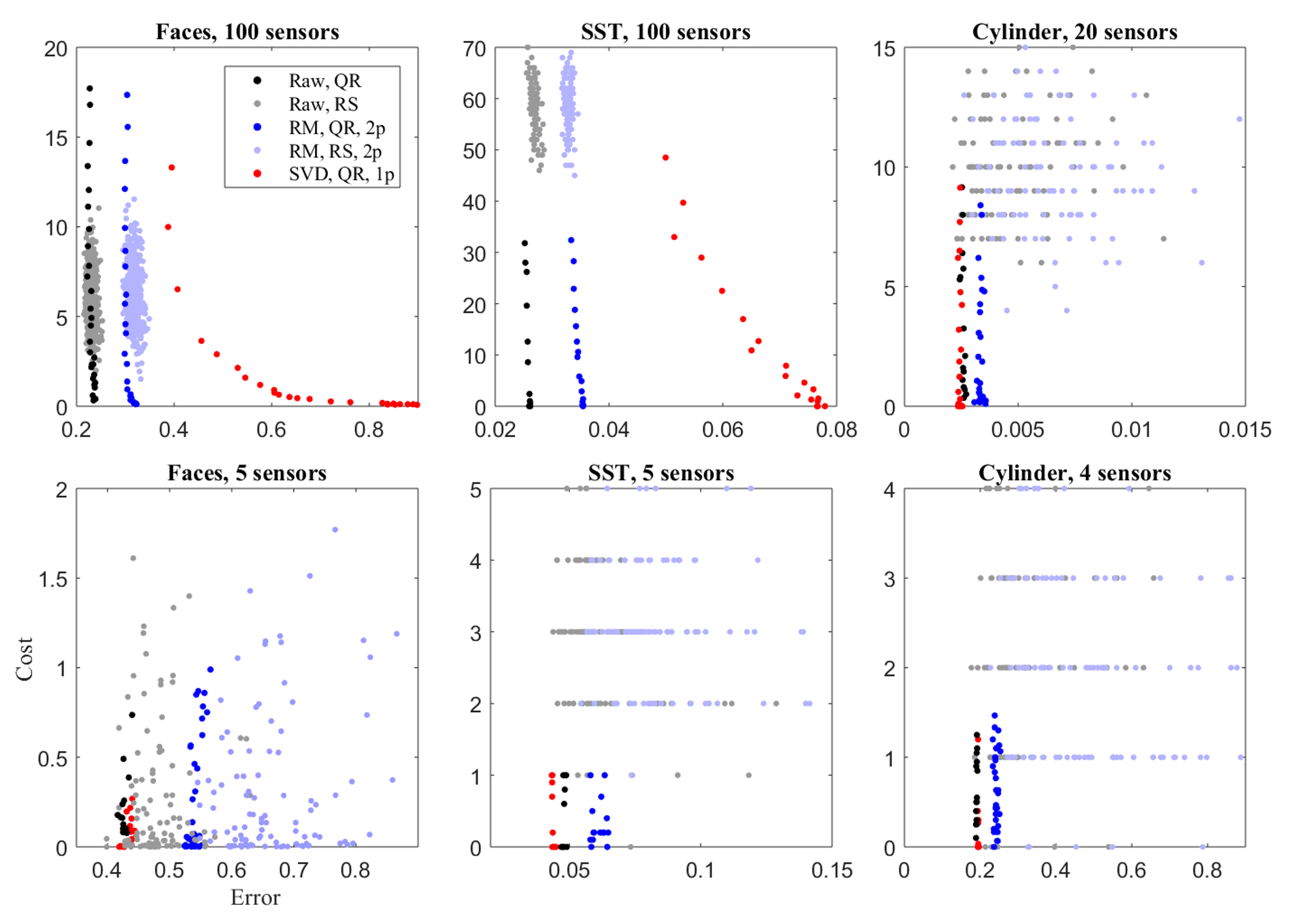}
			\vspace*{-.1in}
	\caption{A comparison of cost versus error results with several different pre-processing methods. The data sets and cost functions tested will be described in more detail in later sections. A Gaussian cost function is used for the eigenface example. The training set is interpolative, not extrapolative, with between 10 and 30 cross validations for each example. The methods used are as follows: First, performing the QR-based sensor placement algorithm directly on the raw data. Second, randomized linear combinations (abbreviated RM in the plot legend, for randomized modes) using $2p$ modes, where $p$ is the number of sensors. Next, QR on $p$ SVD modes. And finally, the sensor placement cost and reconstruction error are calculated from randomly-distributed sensors (abbreviated RS), with either the raw data or $2p$ randomized modes as a reconstruction basis. Randomized sensors using the SVD basis had significantly higher reconstruction errors than any other method, and the results are not shown here. The abbreviations in the legend are explained in Table \ref{table:Comparisons_Colors}.}
	\label{fig:Comparisons}
\end{figure*}

\begin{table}[t]
\begin{center}
\begin{tabular}[h]{|m{2em} | m{4em} | m{15em} |}
\hline
\centering\textbf{Color} &\centering \textbf{Method} &\textbf{Description}\\
\hline\hline
\centering\tikz\draw[black,fill=black] (0,0) circle (.6ex); & \centering Raw, QR & Performing the QR-based algorithm directly on the raw data.\\
\hline
\centering\tikz\draw[mygray,fill = mygray] (0,0) circle (.6ex); &\centering Raw, RS & Randomly-placed sensors, using the raw data as a basis.\\
\hline
\centering\tikz\draw[blue, fill=blue] (0,0) circle (.6ex); &\centering RM, QR, 2p & QR on the randomized modes ${\bf \Psi} = {\bf G}{\bf X}$, ${\bf G} \in \mathbb{R}^{2p\times m}$ a matrix with randomized entries.\\
\hline
\centering\tikz\draw[myblue,fill=myblue] (0,0) circle (.6ex); &\centering RM, RS, 2p & Using randomly-selected sensor locations, and $2p$ randomized modes as a basis.\\
\hline
\centering\tikz\draw[red, fill = red] (0,0) circle (0.6ex); &\centering SVD, QR, 1p & QR on the first $p$ SVD modes.\\
\hline
\end{tabular}
\end{center}
\caption{A brief description of the sensor placement methods used to create Figure ~\ref{fig:Comparisons}, and the colors in which they are plotted.}
\label{table:Comparisons_Colors}
\end{table}

Before proceeding to the cost-constrained placement
examples, we will first briefly discuss the question
of data pre-processing for sensor placement. In the
notation of Section~\ref{sec:prelim}, pre-processing
refers to the process of creating the matrix
${\bf \Psi}^{tr}$ from the training data ${\bf X}^{tr}$ (we then
apply the QR-based algorithm to ${\bf \Psi}^{tr}$). When
selecting $p$ sensors, a common choice for the matrix
${\bf \Psi}^{tr}$ is to set it as the first $p$ right singular
vectors of ${\bf X}^{tr}$
\cite{chaturantabut2010nonlinear,drmac2016new,manohar2017data}.
Inspired by the randomized linear algebra community
\cite{liberty2007randomized,halko2011finding}, we also
consider setting ${\bf \Psi}^{tr} = {\bf G} {\bf X}^{tr}$, where the entries
of ${\bf G} \in \mathbb{R}^{2p\times m}$ are drawn from a standard normal
distribtuion, i.e. we set the rows of ${\bf \Psi}^{tr}$ to be
random linear combinations of the rows of ${\bf X}^{tr}$. We
also consider the performance when setting ${\bf \Psi}^{tr} = {\bf X}^{tr}$,
i.e. the performance without pre-processing. The
number of singular vectors or random linear
combinations used is open-ended; we found that the
choices above gave reasonably optimal performance
for each pre-processing technique.

Our data sets are the Extended
Yale Face Database B, the Optimally Interpolated
Sea-Surface Temperature data set from NOAA, and
simulation data for fluid flow behind a cylinder, all
of which we will describe in more detail in the
next section. For the face and sea surface temperature
data, we consider two types of training sets:
interpolative and extrapolative. By interpolative,
we mean that we have sampled a subset of the data
that draws from all regions of the parameter
space. By extrapolative, we mean that we have purposefully
missed data from a portion of the parameter space.
For the faces data, that means we leave out
all images belonging to 20\% of the individuals. For
the sea surface temperature data, that means we
leave out samples from the last 20\% of the dates.
The data which is left out forms the testing set.
Finding good sensor locations for the extrapolative
training sets is a harder problem; the sensor
locations must reasonably generalize to samples
of data which may be unlike anything in the
training set. We do not make this distinction
for the fluid simulation data as it is much
lower rank and periodic in time.

The top row of Figure~\ref{fig:Interp_Extrap} shows
the spectrum of normalized singular values for all
three data sets. These plots include a gray line at
the Gavish-Donoho \cite{gavish2014optimal} optimal
hard-threshold cutoff, which is an estimate of the
rank beyond which the SVD modes represent additive noise.
It is apparent that the fluid flow data set is
fundamentally different from the other
two, having a sharp elbow at the cutoff, as opposed
to a slow decay.

The remainder of Figure \ref{fig:Interp_Extrap}
plots the relative interpolation error \eqref{eq:errtest}
computed for the test set as a function of the
number of sensors using each of the three methods
for pre-processing (SVD modes, randomized modes,
and the raw data) described above.  
The sea surface temperature
and eigenface data sets both have an interesting
feature not present in the fluid flow data: the error
from the SVD basis has a local minimum at a very small
number of modes -- five for the temperature data and 
approximately ten for the eigenfaces.
While the error for the interpolative training
set begins to slowly decrease again as the number of
sensors is increased, the performance for extrapolative
data gets worse as more sensors are added, up until at
least 150 sensors. This unexpected behavior
reveals that, for systems with slow singular
value decay, there is an overfitting problem when
using SVD modes which occurs well before those
modes correspond to additive sensor noise.

\begin{figure*}[t]
	\centering
	\includegraphics[width = 0.9\textwidth]{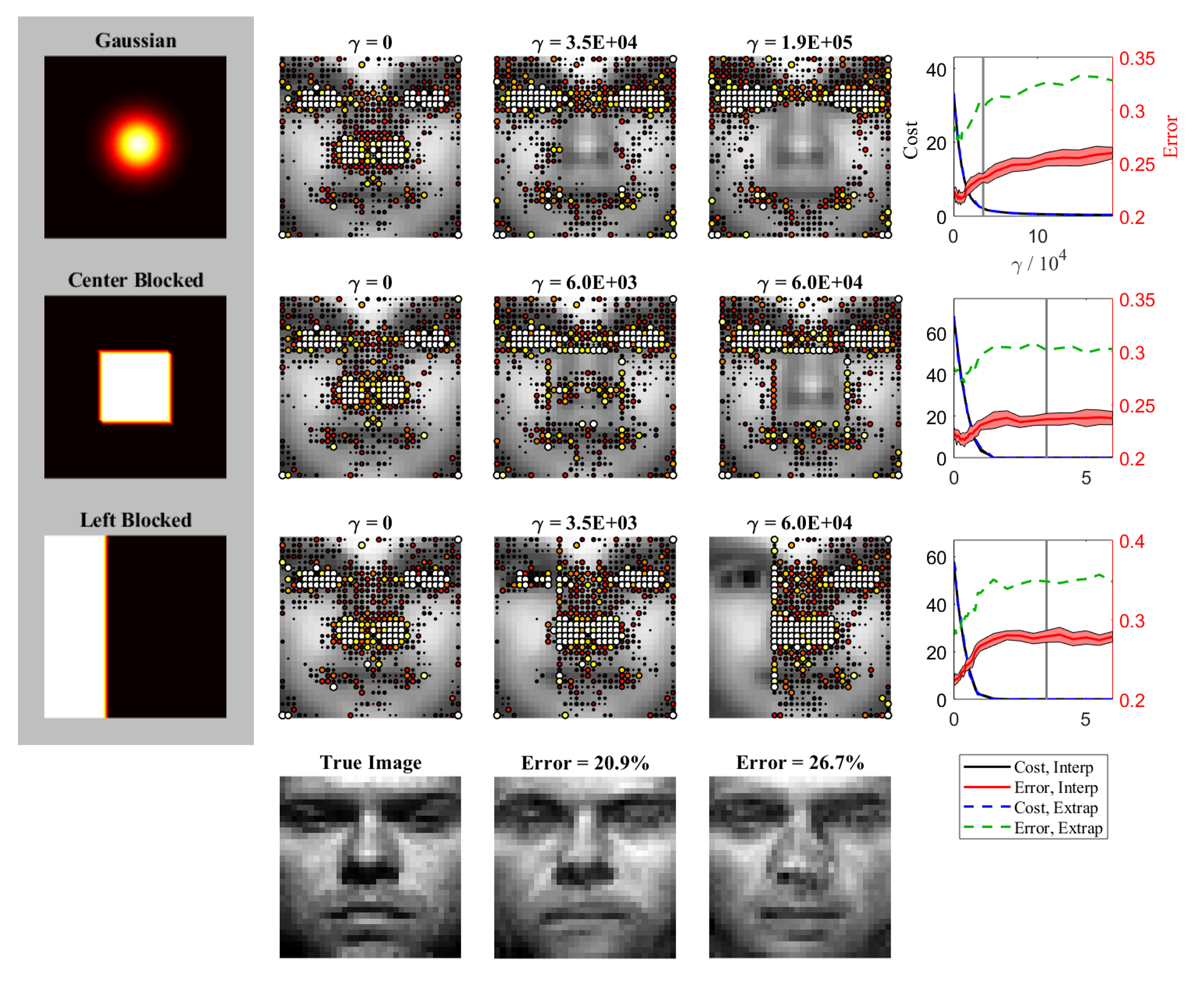}
	\vspace*{-.2in}
	\caption{Average sensor locations for eigenface reconstruction for three different cost functions. The cost functions are plotted in the left 	column, where white indicates regions of highest cost and black shows regions of zero cost. The three central columns show the locations of 200 sensors placed by the QR-based algorithm, averaged over 20 cross validations, for increasing values of the weighting factor $\gamma$. The final column plots cost and error against $\gamma$ for each cost function, with bands indicating error bars. Both interpolation and extrapolation results are given. The vertical gray lines indicate the value of $\gamma$ at which the sensors are plotted in the middle column. The bottom row shows typical example reconstructions of one of the photos, for reference.}
	\label{fig:Face_Sensors}
\end{figure*}

The reconstruction errors for
sensors based on random linear combinations or the
raw data do not have this behavior (except for a weak
effect with extrapolative sea surface temperature data),
nor do any of the cylinder trials.
Indeed, the error for the random linear combinations
and the raw data behaves as expected, decaying at a
rate that follows the error obtained from projecting
the test set onto the first $p$ singular modes of the
training set (this rough indicator of optimal behavior
is plotted as a solid red line in the figure).

We make a few conclusions based on these
pre-processing results.
If the goal of pre-processing is to improve the
quality of the sensors, then it appears that
using no pre-processing is the preferred method,
except when placing a very small number of sensors,
where the SVD mode method displays an advantage.
If the goal of the pre-processing is to
reduce the size of the problem and
speed up the optimization procedure, then it
appears that using randomized linear combinations
of the modes is preferable to using a
limited number of SVD modes. (Note that these
SVD modes would have to be computed with
an accelerated procedure in order
for using them to represent a speed-up over
the QR-based sensor placement algorithm,
again, with the caveat that SVD modes behave
better for a small number of sensors.)

\begin{figure*}[t]
	\centering
	\includegraphics[width = 0.9\textwidth]{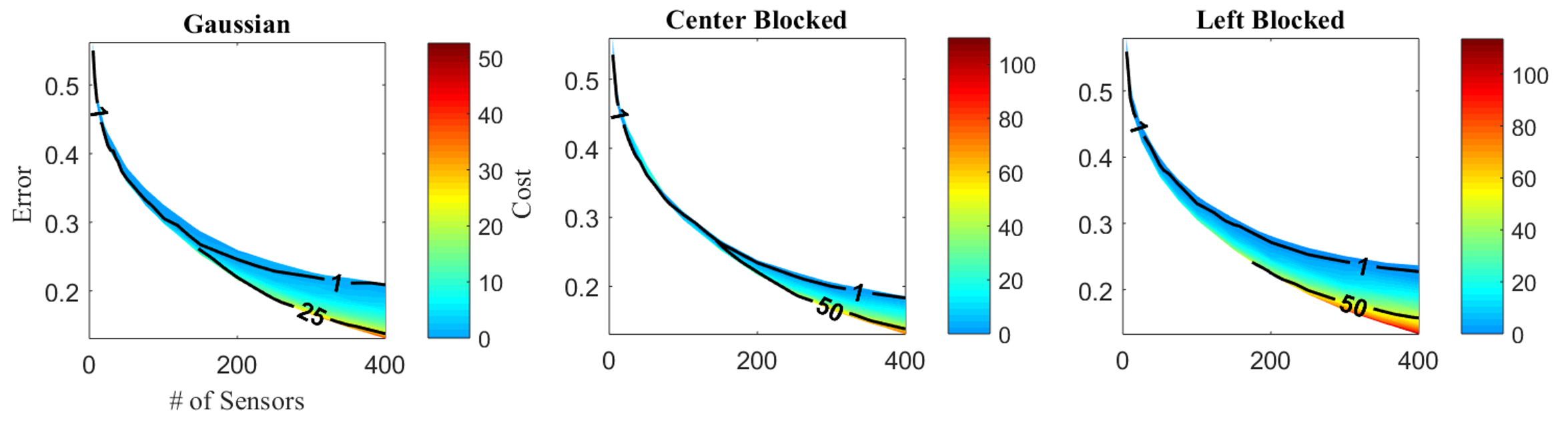}
		\vspace*{-.1in}
	\caption{The cost landscape for eigenface reconstruction for all three cost functions, with cost plotted as a two-dimensional color map against error and the number of sensors. Contours indicate lines of constant cost.}
	\label{fig:Face_Contour}
\end{figure*}

\section{Applications} \label{sec:results}

When factoring in the effect of cost, we observe similar
behavior to the cost-free case analyzed in the previous
section. Figure \ref{fig:Comparisons} provides an overview
of the performance for the various pre-processing methods,
now with a non-zero cost associated to each location. In the
plots, we generate several cost, error pairs for each
pre-processing technique by varying the cost function
weighting $\gamma$ in the QR-based algorithm. We also
plot the performance of randomly drawn sensors for the
sake of comparison. All three data sets are considered,
with both a large and small number of sensors. See
Table \ref{table:Comparisons_Colors} for details on the
figure labels.

For both the eigenface and sea surface
temperature data sets, using the raw data at 100 sensors leads
to the lowest error at a comparable cost, and randomized linear 
combinations with $2p$ modes gives the next lowest error,
followed by SVD with $p$ modes. At 5 sensors, the latter
is reversed, with SVD performing comparably to or better than the raw data,
as is the case with both trials for the fluid flow behind a cylinder.
We observe that our QR-based method outperforms the
best randomly chosen sensors when placing a large
number at low cost, while choosing the best random
sensors may offer an advantage when placing a small
number (the randomly chosen sensors are more likely
to contain the brute force answer in this case).
We note that the randomized data performs worse than
the raw data in all examples, but the behavior is
consistent and may be worth the reduced computational
cost in some applications.

In the remainder of this section, we more closely
analyze the performance of the algorithm for each
of our data sets. For brevity, we consistently
use a randomized linear combination of data vectors for the
pre-processing technique in these examples, noting
that the behavior described above would be maintained
if all techniques were tested.

\subsection{Eigenfaces}

\begin{figure}[t]
	\centering
	\includegraphics[width = 0.9\columnwidth]{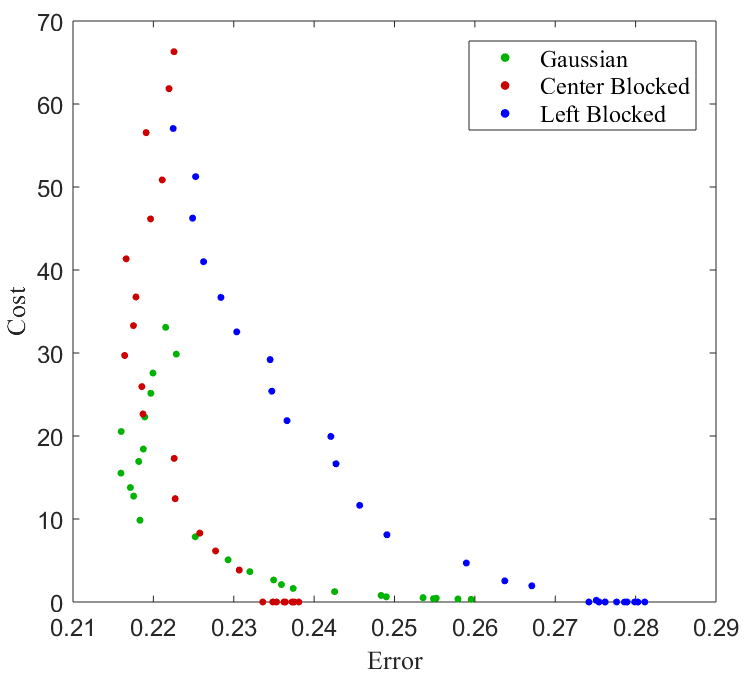}
	\vspace*{-.1in}
	\caption{Cost versus error for eigenface reconstruction with three different cost functions. Results are given for the case of 200 sensors.}
	\label{fig:Face_CvE_200}
\end{figure}

\begin{figure}[t]
	\centering
	\includegraphics[width = 0.9\columnwidth]{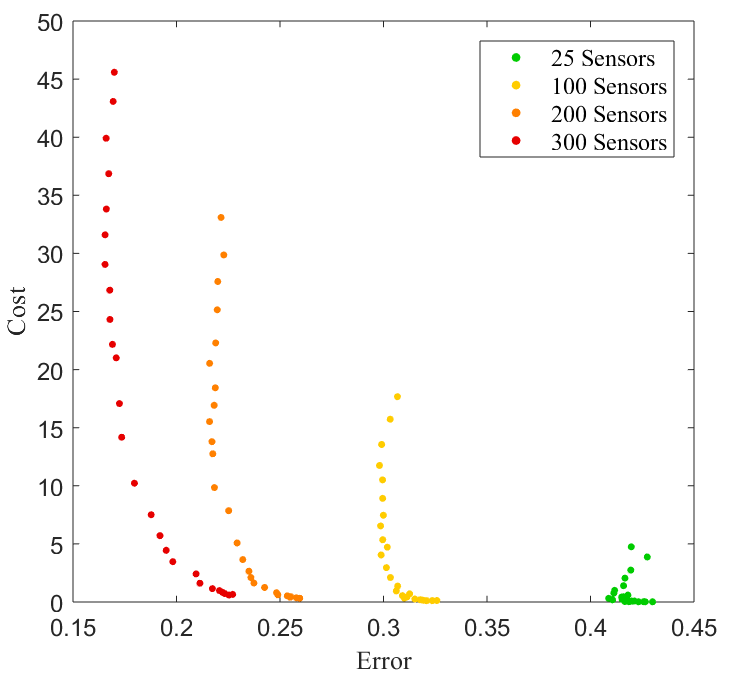}
	\vspace*{-.1in}
	\caption{Cost versus error for eigenface reconstruction with a Gaussian cost function, for several different numbers of sensors.}
	\label{fig:Face_CvE_Gauss}
\end{figure}

\begin{figure*}[t]
	\centering
	\includegraphics[width = 0.9\textwidth]{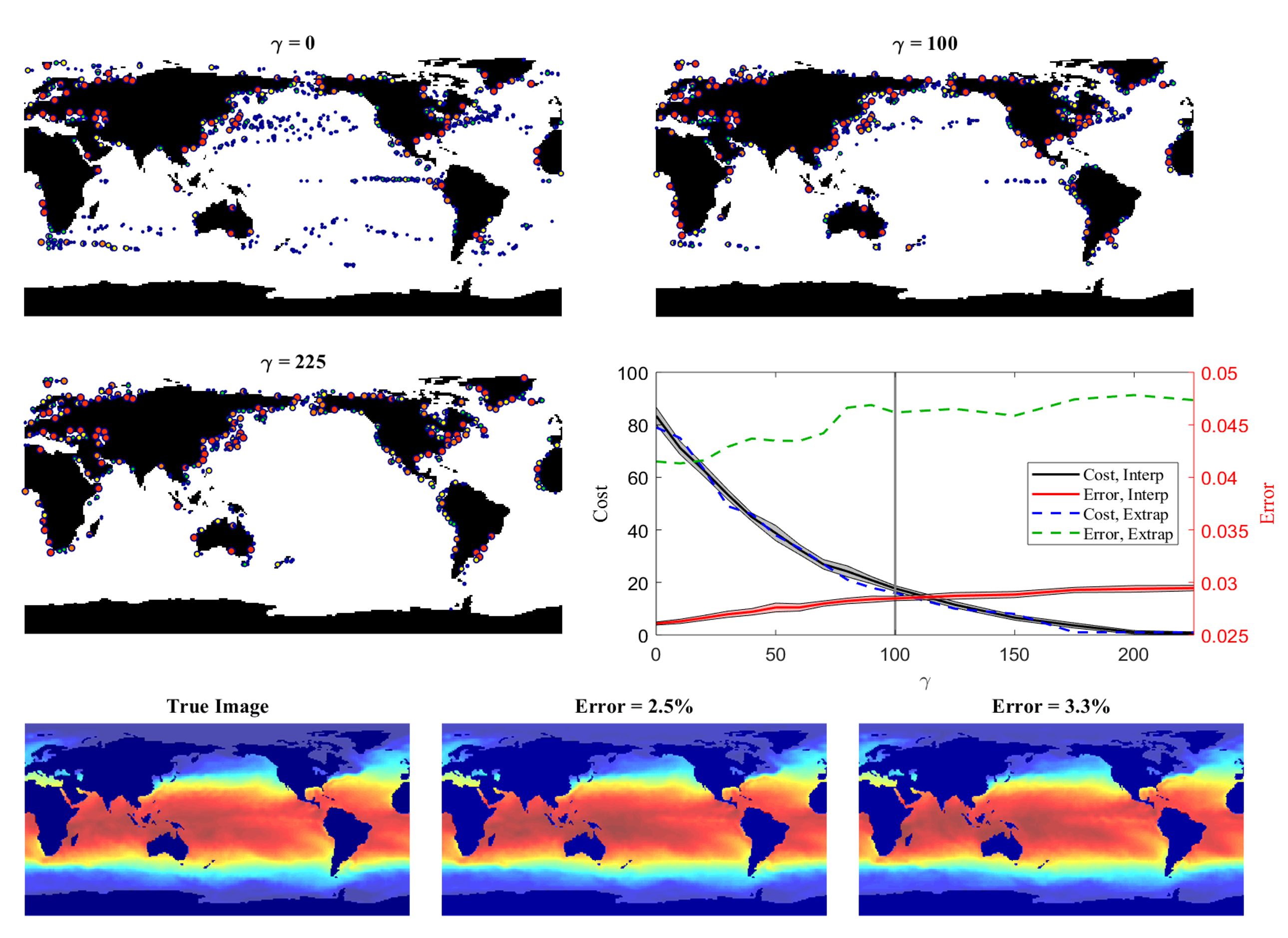}
	\vspace*{-.1in}
	\caption{Sensor locations for sea surface temperature reconstruction with 200 sensors. The cost function considered was a step function which was zero up to two pixels off land and equal to one everywhere else. Locations are shown for three different values of the weighting factor $\gamma$, and are averaged over ten cross validations. Size and color of a data point indicate the frequency with which a sensor was placed there. The fourth image plots cost and error against $\gamma$, for both interpolative and extrapolative data sets. The bottom row shows a comparison of an example temperature snapshot along with two reconstructions of it yielding two different accuracies.}
	\label{fig:SST_Sensors}
\end{figure*}

\begin{figure}[t]
	\centering
	\includegraphics[width = 0.9\columnwidth]{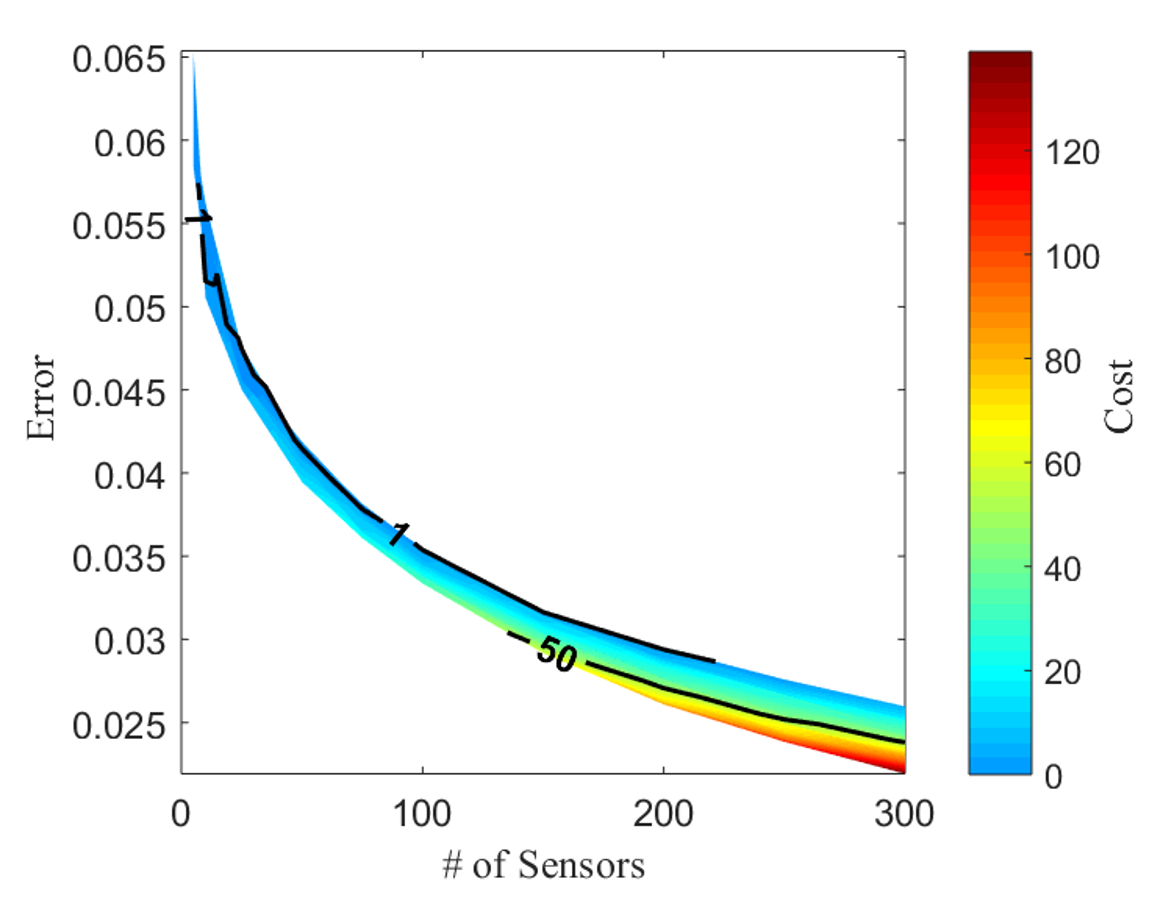}
	\vspace*{-.1in}
	\caption{The cost landscape for sea surface temperature reconstruction, plotted as a color map against error and the number of sensors. Contours show lines of constant cost.}
	\label{fig:SST_Contour}
\end{figure}

The algorithms are first tested on the Extended Yale Face Database B, referred to here as the eigenface data set \cite{CHHH07,CHH07b,CHHZ06,HYHNZ05}. It comprises about 64 images each of 38 individuals under various lighting conditions. The images are downsized to $32\times32$ pixels.

Unless otherwise stated, the tests are conducted on an interpolative training set, by randomly selecting 80\% of the images. The cost weighting factor $\gamma$ is then progressively increased as the sensors are placed using the QR-based algorithm. Once the sensor locations are selected, the cost of the given array is calculated, and their average reconstruction accuracy is evaluated using the remaining 20\% of the photographs. A 20-fold cross validation is performed.

Three cost functions are tested: 1) a Gaussian function, such that it is most expensive to place sensors in the center of the face, 2) a step function uniformly penalizing sensors in the middle ninth of the photographs, and 3) a step function penalizing the left third of the data set.

\begin{figure*}
	\centering
	\includegraphics[width = 0.9\textwidth]{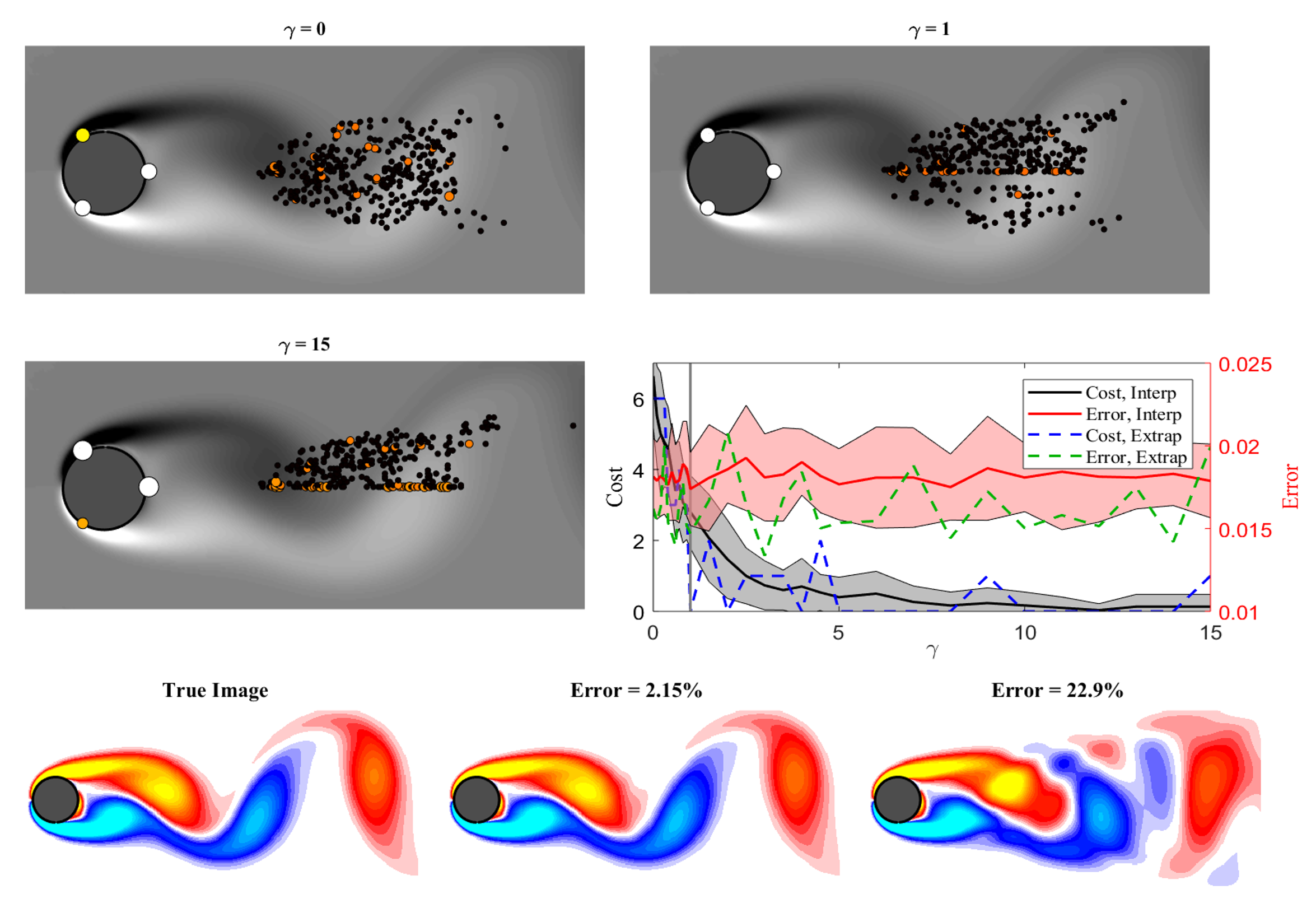}
	\vspace*{-.1in}
	\caption{Sensor locations and reconstructions of the flow behind a cylinder, obtained by the QR-based algorithm, using 14 sensors. The cost function was a step function blocking the lower half of the domain. The first three plots show the placements of the sensors averaged over 30 trials. The fourth shows cost and error plotted against $\gamma$, where the bands indicate error bars. The bottom row gives a comparison of the true image of a snapshot along with two different reconstructions.}
	\label{fig:Cyl_Sensors}
\end{figure*}

All three cost functions and a few corresponding sensor arrays are shown in Figure \ref{fig:Face_Sensors}. These are the average locations over the twenty cross-validation runs, generated by the QR algorithm with 200 sensors, shown as a scatter plot on top of the dominant eigenface mode. Marker size and color indicate the frequency with which a sensor was placed at a given location. As expected, when the cost function weighting is increased, sensors are gradually pushed out of the regions of higher cost. This allows the total cost to be lowered at the expense of decreasing reconstruction accuracy, as demonstrated in the right-hand column of the figure, which plots cost and error on separate axes, as a function of $\gamma$. Note that $\gamma$ is an arbitrary weighting, and the same value of $\gamma$ can have very different effects depending on the cost function. Extrapolative cost and error are shown in the same panel, where the extrapolative error is higher than the interpolative, at an identical cost. Additionally, the bottom row of the figure shows several reconstruction examples for one of the faces.

For many practical applications of these methods, the final goal will be to minimize reconstruction error given a predetermined budget. To that end, cost-error landscapes are constructed by calculating sensor array cost and reconstruction error for different numbers of sensors. The results are shown in Figure \ref{fig:Face_Contour}, which shows the landscapes for each cost function. These landscapes are plotted as a color map according to cost. Cost contours on this color map represent a hypothetical budget, so the optimum configuration for a given budget can be found by following the appropriate contour to the lowest possible error. Note that the upper edges of the contour plots indicate the minimum cost and maximum error for a given number of sensors, and the lower edges indicate the minimum error and maximum cost.

Cross sections of the cost landscapes as plots of cost versus error are given in Figures \ref{fig:Face_CvE_200} and \ref{fig:Face_CvE_Gauss}. The former shows such cross sections for each of the three cost functions, using 200 sensors, and the latter shows cost versus error for a Gaussian cost function with 25, 100, 200, and 300 sensors.

\subsection{Sea surface temperature}

The next data set we consider is the NOAA\_OISST\_V2 mean sea surface temperature set \cite{noaa2018, banzon2016, reynolds2007}, comprising weekly global sea surface temperature measurements between the years of 1990 and 2016. There are a total of 1400 snapshots on a $360\times180$ spatial grid. The QR-based algorithm is trained on 1100 randomly-selected snapshots and tested on the remaining 300. Ten cross validations are performed. The cost function used is a step function which penalizes placing sensors too far from shore, being zero for locations one and two pixels off land, and equal to one everywhere else.

Average sensor locations over ten cross validations with 200 sensors are shown in Figure \ref{fig:SST_Sensors}, as a scatter plot where the size and color of a data point indicate the frequency with which a sensor was placed in that location (blue being the least frequent, red being the most frequent). As the cost function weighting is increased, more sensors move within the unblocked regions offshore, until the cost reaches zero. Plots of cost and error as functions of $\gamma$ are given in the fourth panel, and the bottom row shows two example reconstructions. As with the eigenfaces, the interpolative trial has a much lower error than the extrapolative trial. Notice that although the reconstruction error increases as cost decreases, here the error does not even reach 3\%, even when the cost is zero.

\begin{figure}[t]
	\centering
	\includegraphics[width =  0.9\columnwidth]{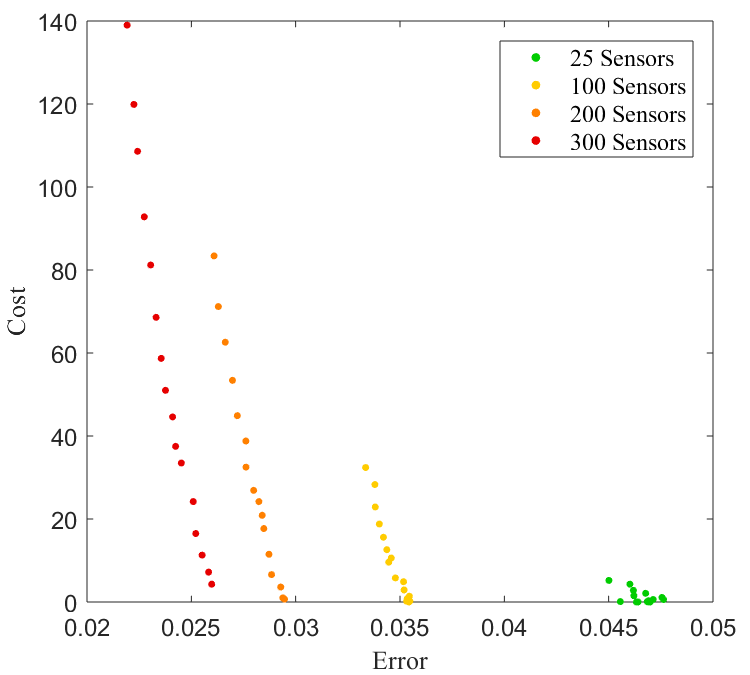}
	\caption{Plots of cost versus error for sea surface temperature reconstruction, with varying numbers of sensors.}
	\label{fig:SST_CvE}
\end{figure}

The landscape of cost as a function of error and the number of sensors is shown in Figure \ref{fig:SST_Contour}, again visualized as a color map with contours of constant cost.
Cross sections of cost versus error for 25, 100, 200, and 300 sensors are shown in Figure \ref{fig:SST_CvE}.

\subsection{Fluid flow around a cylinder}

As a final example, vortex shedding of a fluid flowing around a stationary cylinder is considered. This data set is low-rank, periodic, and vertically symmetric, making it a significant contrast to the previous two examples. The flow data were generated using the immersed boundary projection method~\cite{taira:07ibfs,taira:fastIBPM} to numerically simulate the Navier-Stokes equations with Reynolds number 100. There are 151 snapshots in time, each on a $199\times449$ spatial grid.

\begin{figure}[t]
	\centering
	\includegraphics[width = 0.9\columnwidth]{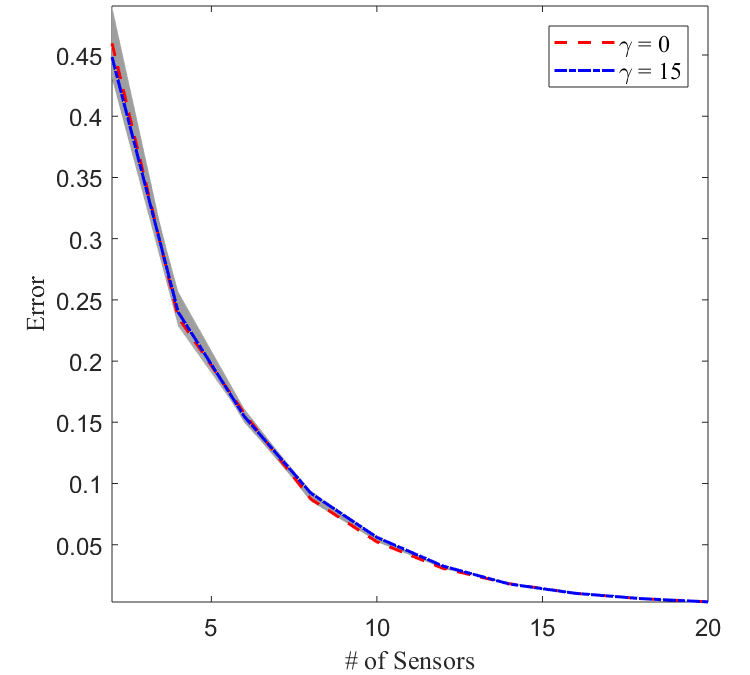}
	\caption{The cost landscape for reconstruction of a fluid flow behind a cylinder. Instead of the color maps made for the previous two data sets, curves of error versus the number of sensors are shown for high and low values of $\gamma$. The gray band indicates the maximum variation in the error.}
	\label{fig:Cyl_Contour}
\end{figure}

\begin{figure}[t]
	\centering
	\includegraphics[width = 0.9\columnwidth]{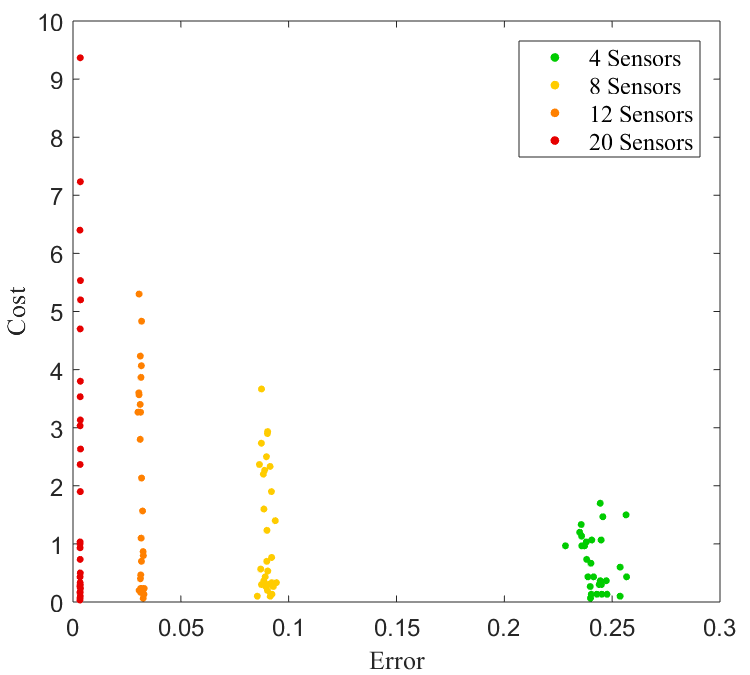}
	\caption{Plots of cost versus error for cylinder flow reconstruction, with varying numbers of sensors.}
	\label{fig:Cyl_CvE}
\end{figure}

The QR-based method is applied to the data set with a cost function that is uniformly one in the lower half of the domain and zero in the upper half. This allows the algorithm to take full advantage of the symmetry of the fluid flow, as can be seen in the figures. The data was generated using an interpolative training set of 120 randomly-selected snapshots; 30 cross validations were performed.

In Fig.~\ref{fig:Cyl_Sensors}, sensor locations for several values of the cost function weighting $\gamma$ are shown. The locations were picked using the QR-based algorithm with 14 sensors, and are averaged over the 30 cross validations, then graphed as a scatter plot on top of an example fluid flow snapshot. The size and color of a data point indicate how frequently a sensor was placed at its location. As expected, when $\gamma$ is increased, the sensors migrate until they occupy the upper half of the plane.

The fourth panel of the figure plots the cost and error on separate axes, as functions of $\gamma$. Because of the symmetry of the data set, the reconstruction error is essentially unchanged with $\gamma$, even as the cost goes to zero. Furthermore, because the flow is periodic, the extrapolative data performs slightly better than the interpolative data. The figure's bottom row shows example reconstructions of a snapshot.

The cost-error landscape is shown in Fig.~\ref{fig:Cyl_Contour}. Notice that there is hardly any variation in the error, so rather than plotting cost as a color map, the error versus number of sensor curves for $\gamma=0$ (highest cost) and $\gamma=15$ (lowest cost) are shown. These curves are essentially identical, further emphasizing that for this particular data set and cost function, the total cost can be lowered with no penalty to the error.
Cost versus error plots are given in Fig.~\ref{fig:Cyl_CvE}, which shows results for 4, 8, 12, and 20 sensors. Regardless of the number of sensors, these cross sections are essentially vertical lines, within the standard deviation of the error.

\section{Conclusion and future directions}
\label{sec:conclusion}

We have developed a QR-based greedy algorithm to place sensors for reconstruction with a cost constraint on sensor locations. This algorithm is tested on three different data sets, eigenfaces, weekly sea surface temperature data, and vortex shedding of a fluid flowing around a cylinder. In all cases, the method is demonstrated to be capable of lowering sensor cost at the expense of marginal increases in reconstruction error.

It is also shown that with or without the inclusion of a cost function, data sets with slow singular value decay have better results by pre-processing the data through a randomized linear combination of modes, rather than through SVD-based rank reduction. Random linear combinations lead to significantly lower reconstruction errors, except at a very low number of sensors.

In fact, for these data sets with slow singular value decay, SVD modes behave in an unexpected way at a low number of modes. The reconstruction error decreases sharply, even surpassing the error obtained by using the full raw data in the case of sea surface temperatures, before increasing again as more sensors are added. This suggests that the SVD is overfitting well before the Gavish-Donoho cutoff, an idea which warrants further exploration in future work. The results also imply that there may be some other pre-processing method which can take advantage of both the SVD behavior at a low number of sensors and the random linear combination or raw data behavior at a higher number of sensors. This will also be explored in the future.

In the meantime, the algorithm presented here provides a way to place sensors under a cost constraint, which could have applications in manufacturing, atmospheric sensing, fluid flow sensing, and many more fields.  Specifically, the algorithm allows one to address three critical engineering design principles in regard to sensors placement:  (i)  For a fixed budget of sensors, where are the best measurement locations, (ii)  What is the minimal number of sensors required to achieve a given reconstruction error, and (iii) How well can inaccessible regions be reconstructed in practice.  Depending upon the application, one or all of these questions may be of central concern.  The computationally tractable approach presented here provides a principled mathematical method for answering these questions.

\section*{Acknowledgment}

T.~Askham would like to thank Mark Tygert for
useful discussions of the interpolative 
decomposition.  
T.~Askham and J. N. Kutz acknowledge support from the Air Force Office of Scientific Research (FA9550-15-1-0385).
S. L. Brunton acknowledges support from the Air Force Research Laboratory (FA8651-16-1-0003) and the Air Force Office of Scientific Research through the Young Investigator Program (FA9550-18-1-0200).

\ifCLASSOPTIONcaptionsoff
  \newpage
\fi



\bibliographystyle{IEEEtran}
\bibliography{IEEEabrv,refs}
%



%

\begin{IEEEbiography}[{\includegraphics[width=1in]{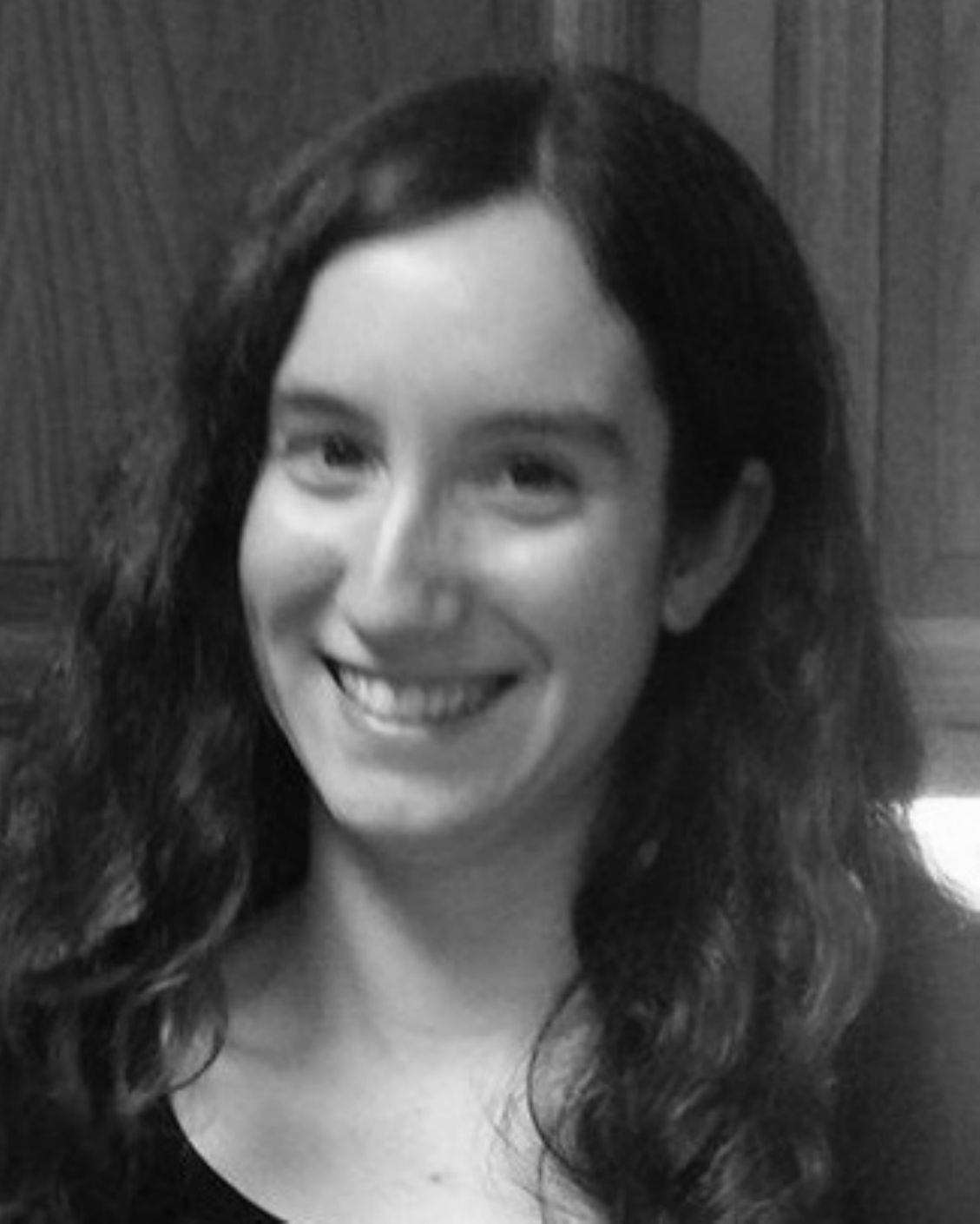}}]{Emily Clark}
received the B.S. degree in physics from Bates College, Lewiston, ME, in 2015. She is a Ph.D. student of physics at the University of Washington.
\end{IEEEbiography}

\begin{IEEEbiography}[{\includegraphics[width=1in]{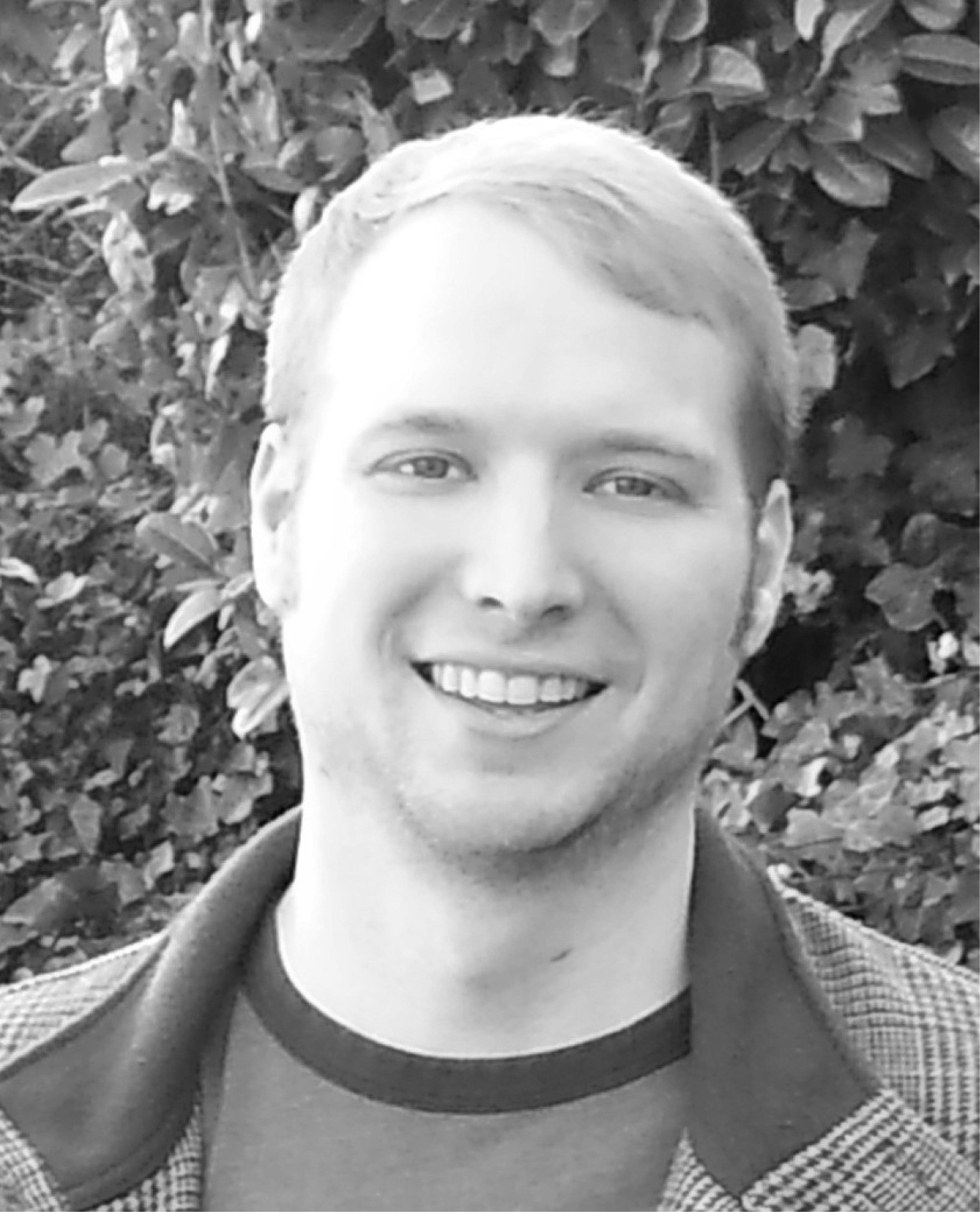}}]{Travis Askham}
received the B.S. degree in applied mathematics and the M.A.
degree in mathematics from the University of California,
Los Angeles, Los Angeles, CA, in 2010 and the Ph.D. degree in 
mathematics from the Courant Institute of Mathematical Sciences 
at New York University, New York, NY, in 2016. He is a 
Research Associate of applied mathematics at the University
of Washington.
\end{IEEEbiography}

\begin{IEEEbiography}[{\includegraphics[width=1in]{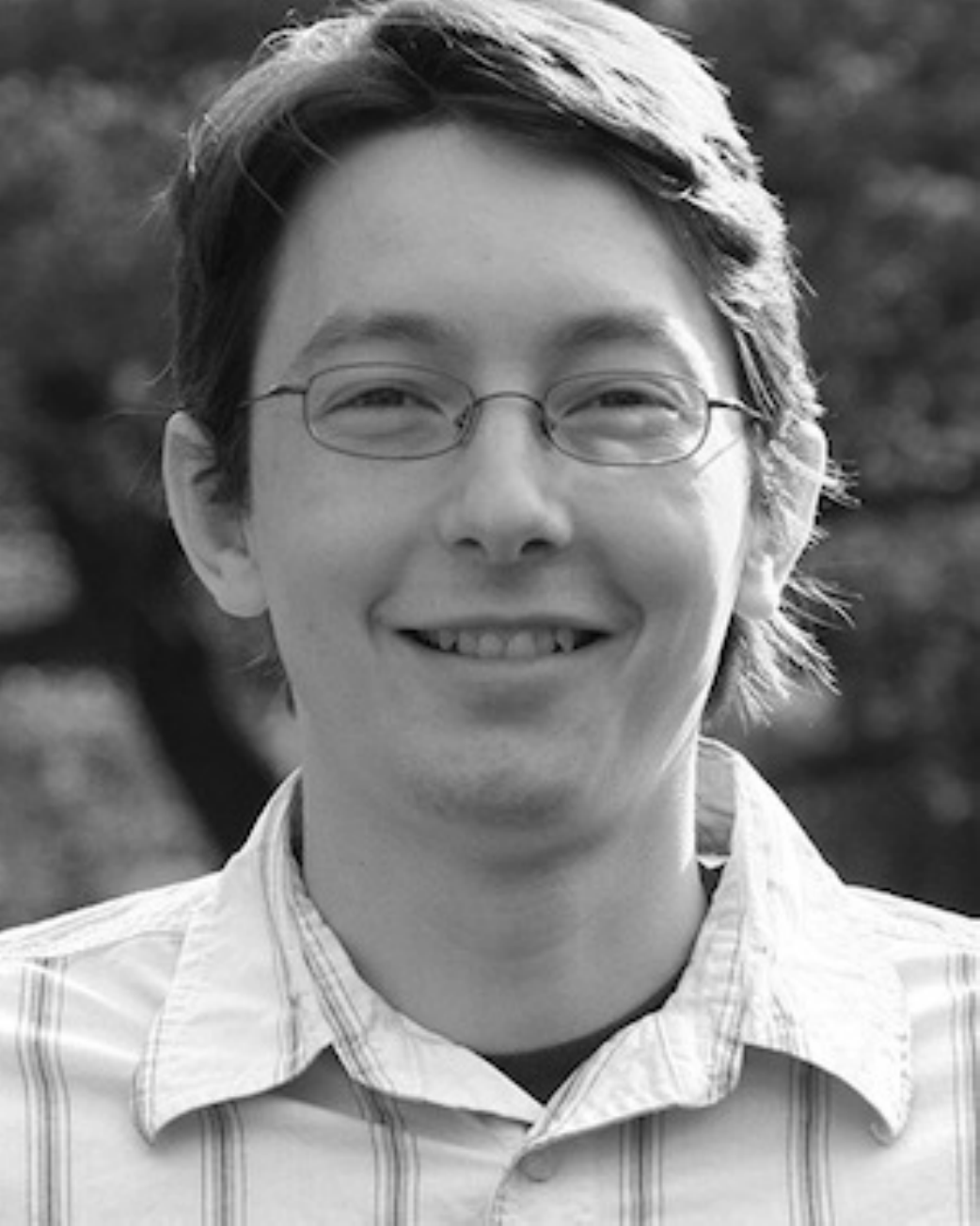}}]{Steven~L.~Brunton}
received the B.S. degree in mathematics with a minor 
in control and dynamical systems from the California 
Institute of Technology, Pasadena, CA, in 2006, and 
the Ph.D. degree in mechanical and aerospace engineering 
from Princeton University, Princeton, NJ, in 2012.  He 
is an Assistant Professor of mechanical engineering and 
a data science fellow with the eScience institute at the 
University of Washington.
\end{IEEEbiography}

\begin{IEEEbiography}[{\includegraphics[width=1in]{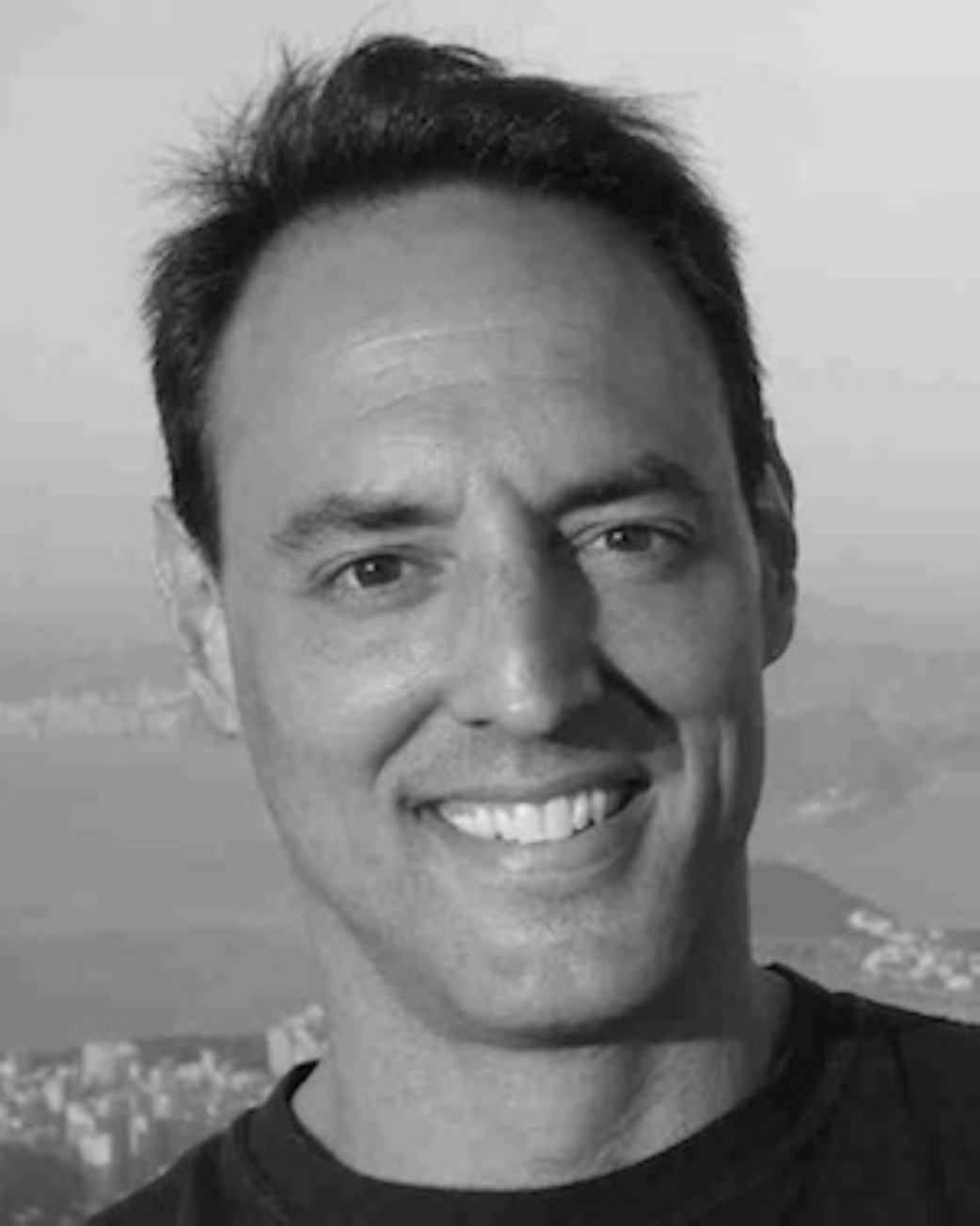}}]{J.~Nathan Kutz}
received the B.S. degrees in physics and mathematics 
from the University of Washington, Seattle, WA, in 1990, 
and the Ph.D. degree in applied mathematics from 
Northwestern University, Evanston, IL, in 1994.  He is 
currently a Professor of applied mathematics, adjunct 
professor of physics and electrical engineering, and a 
senior data science fellow with the eScience institute 
at the University of Washington.
\end{IEEEbiography}







\end{document}